# ORGANIZATION OF PARAMETER SPACE FOR SIMPLE CIRCLE MAPS: THE FAREY WEB

K. BRUCKS, J. RINGLAND, AND C. TRESSER

ABSTRACT. We define the Farey web — a collection of loci in the parameter plane of families of simple non-invertible maps of the circle. We prove some results about the arrangement of these loci and their relationships with other dynamically significant features of the parameter plane. The results enable us to provide short proofs for a number of theorems about the organization of frequency-locking.

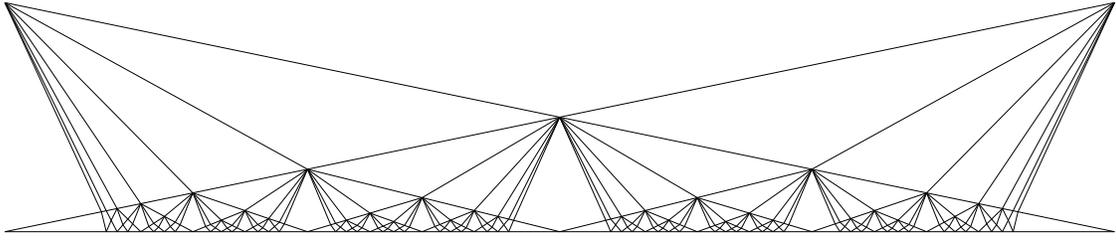

**Figure 1.** The Farey web under construction.

## 1. INTRODUCTION

In this note we define the *Farey web*, a collection of loci in the parameter plane of certain nice families of non-invertible maps of the circle. We investigate the arrangement of the strands of the web, mutually and in relation to other dynamically significant loci. The strands are the loci of existence of certain orbits, and the junctions of the strands are special points of the frequency-locking regions (regions where all orbits rotate at the same average rate). Using our results about the web, we describe a geometrical construction (see Figure 1) of a homeomorphic image of its simplest possible realization in the class of maps we consider, and we prove very simply several theorems about the organization of frequency-locking. Specifically, we generalize a result of [BJ, FT] about the lack of local connectedness of the boundary of zero-width rotation interval, and supply proofs which have simple graphical interpretations of two results of [BT], one strengthened slightly, about the relationships between frequency-locking regions of different frequencies. It is our hope that the Farey web will yet further aid in the description of the transition from simple to complicated dynamics in these two-parameter families of circle maps.

In Section 2 we state our results after first providing the definitions necessary to do so. Section 3 summarizes some knowledge about degree-one lifts in preparation for Section 4 which contains our proofs.

---

## 2. Summary of Results

Let $F$ be a *degree-one lift*, i.e., $F : \mathbb{R} \to \mathbb{R}$ is a continuous map such that $F(x + 1) = F(x) + 1$ for all $x \in \mathbb{R}$. Hence, $F$ is a lift to the universal cover $\mathbb{R}$ of some degree-one endomorphism $f$ of the circle $\mathbf{T} = \mathbb{R}/\mathbb{Z}$, i.e., $\pi \circ F = f \circ \pi$, where $\pi : \mathbb{R} \to \mathbf{T}$ is given by $\pi(x) = x \bmod 1$. An important characteristic of a lift $F$ is its *rotation set*, $Rot(F)$. To define it, for each $x \in \mathbb{R}$ we set

$$\underline{\rho}_F(x) = \liminf_{m \to \infty} \frac{F^m(x) - x}{m} \quad \text{and} \quad \overline{\rho}_F(x) = \limsup_{m \to \infty} \frac{F^m(x) - x}{m},$$

where as usual $F^0 = Id$ and $F^m = F \circ F^{m-1}$. Then we define,

$$\underline{Rot}(F) = \{\underline{\rho}_F(x) \mid x \in \mathbb{R}\},$$
$$\overline{Rot}(F) = \{\overline{\rho}_F(x) \mid x \in \mathbb{R}\},$$
$$Rot(F) = \{\rho_F(x) \mid x \in \mathbb{R} \text{ and } \rho_F(x) = \underline{\rho}_F(x) = \overline{\rho}_F(x)\}.$$

It is known that $\underline{Rot}(F) = \overline{Rot}(F) = Rot(F)$ and that this set is either a single point or a non-degenerate closed interval (see *e.g.* [ALM, Theorem 3.7.20]).

*Remark.* Throughout this note we make references to [ALM] for the convenience of the reader and the authors; the reader will find original references traced and discussed in [ALM].

The map $F$ is said to be *frequency-locked* if $Rot(F)$ is a single point. Any degree-one lift which is non-decreasing is frequency-locked; this is a simple extension of the classical theory of rotation numbers developed by Poincaré [P]. For a degree-one lift $F$, its *upper* and *lower monotone bounds*, $F_l$ and $F_u$, are defined as follows (see Figure 2(a)):

$$F_l(x) = inf\{F(y) \mid y \geq x\} \quad \text{and} \quad F_u(x) = sup\{F(y) \mid y \leq x\}.$$

These are useful because $Rot(F) = [\rho_{F_l}, \rho_{F_u}]$ [ALM].

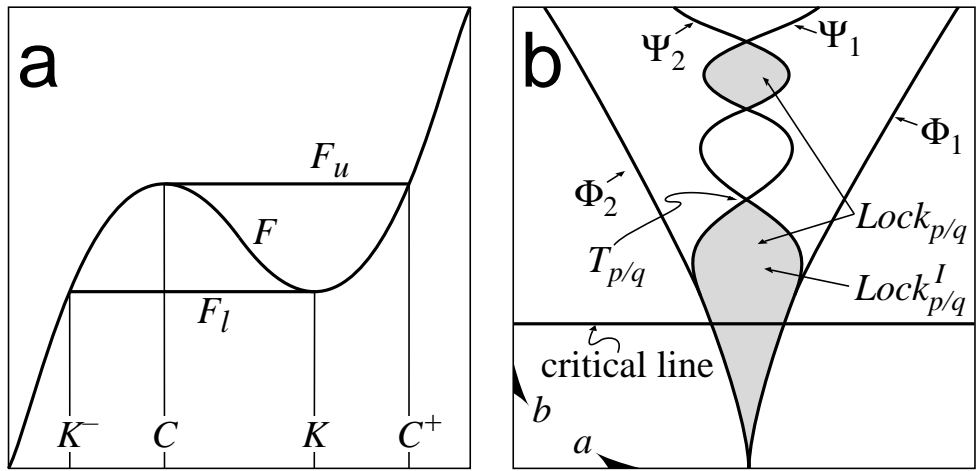

**Figure 2.**

As in [B, FT, MT] we consider degree-one circle maps which are $C^3$, piecewise strictly monotone, have at most two critical points, both turning points when there are two, and have negative Schwarzian derivative (where defined) when they do have two critical points. Recall that the *Schwarzian* derivative of $F$ is defined as:

$$S(F) = \frac{F'''}{F'} - \frac{3}{2}\left(\frac{F''}{F'}\right)^2.$$

We do not use the negativity of this quantity directly, but rather through the implications described in Facts 1(v), 1(vi), 6 and 7. We say a lift $F$ of a such a circle map belongs to Class $\mathcal{B}$. When $F$ has two turning points per fundamental domain, we then choose the origin such that the locations of the maximum and the minimum in $[0, 1)$, denoted by $C$ and $K$ respectively, satisfy $C < K$. Additionally we denote by $K^-$ the closest point to $K$ below $K$ such that $F(K^-) = F(K)$, and by $C^+$ the closest point to $C$ above $C$ such that $F(C^+) = F(C)$. See Figure 2(a).

We restrict our attention to two-parameter families $\{f_{a,b}\}$ of circle maps such that a corresponding family of lifts $\{F_{a,b}\}$ satisfies the following conditions.

(i) Each lift $F_{a,b}$ is in class $\mathcal{B}$.
(ii) Both $F_{a,b}$ and $(d/dx)F_{a,b}$ depend continuously on $(a, b)$.
(iii) The parameter $a$ is a translation parameter as, for example, in the *standard sine family*, $F_{a,b}(x) = x + a + \frac{b}{2\pi} \cdot \sin(2\pi x)$. Actually, [FT] only require that $a_2 > a_1$ implies that $F_{a_2,b} > F_{a_1,b}$.
(iv) The (Hausdorff) distance, $\Delta$, between the upper and lower monotone bounds is a nondecreasing function of $b$ such that the value $0$ can be attained and the value $1$ can be exceeded. Note that no frequency-locking is possible for $\Delta > 1$.

As a consequence of the above conditions, there exists a well-defined *critical line* in parameter space, $b = max\{b|\Delta = 0\}$, ($b = 1$ in the standard sine family, for example). On this line, $F_{a,b}$ has a unique critical point in $[0, 1)$, and above it, it has two of them. On the critical line we set $K^- = C = K = C^+$. We say $\{f_{a,b}\}$, or $\{F_{a,b}\}$, belongs to Class $\mathcal{N}$ (for nice) if conditions (i) through (iv) are satisfied. We assume that in a family $\{F_{a,b}\}$ the origin is chosen so that when $0 \leq \Delta < 1$, we have $0 < K^- \leq C \leq K \leq C^+ < 1$. Note that although $K^-$, $C$, $K$, and $C^+$ do not depend on the parameter $a$, they do depend on the parameter $b$; but to keep notation simple, we do not usually use a subscript to indicate this.

For a two-parameter family of degree-one lifts $\{F_{a,b}\}$ in Class $\mathcal{N}$ and $\omega \in \mathbb{R}$ the *Arnold tongue*, $A_\omega$, is defined as the subset of parameter space where $\omega$ belongs to the rotation set. Due to condition (iv) above, for each $\omega$, $A_\omega$ is a connected region. The *frequency-locking region* with frequency $\omega$ is the subset $Lock_\omega \subset A_\omega$ of parameter space where the rotation set is $\{\omega\}$, and the height of $Lock_\omega$ is the largest value of $b$ such that some map $F_{(a,b)}$ belongs to $Lock_\omega$. (Such a value of $b$ exists as a consequence of condition (iv) above.) It has recently been shown [EKT] that for the standard sine family the frequency-locking region for each frequency is connected. However, in class $\mathcal{N}$ one can construct families that exhibit "pathologies" such as are sketched in Figure 2(b). Hence, we denote by $Lock_\omega^I$ the connected component of $Lock_\omega$ that contains all invertible maps in $Lock_\omega$. The *reduced height* of $Lock_\omega$ is used to mean the height of $Lock_\omega^I$, i.e, the largest value of $b$ such that some map $F_{(a,b)}$ belongs to $Lock_\omega^I$. For $\omega \in \mathbb{R}/\mathbb{Q}$, the reduced height of $Lock_\omega$ is the critical line value. By a *tip* of the frequency-locking region $Lock_{\frac{p}{q}}$ ($\frac{p}{q} \in \mathbb{Q}$), we mean a point $(a_0, b_0)$ such that $(a, b_0) \in Lock_{\frac{p}{q}}$ iff $a = a_0$. The greatest value of b in $Lock_{\frac{p}{q}}^I$ is achieved at a tip, and we denote this point by $T_{\frac{p}{q}}^I$. See Figure 2(b).

**Fact 1** [B, Lemmas 2.5, 5.1, proof Theorem 3.1]. Let $F$ be a degree one lift such that $F$ is strictly piecewise monotone with precisely two turning points in $(0,1)$ (a maximum at $K$ and minimum at $C$). For $a \in \mathbb{R}$, let $F_a = F + a$. Fix $\frac{p}{q} \in [0,1]$. Then there exist closed intervals $[\Psi_1, \Phi_1] \subset \mathbb{R}$ and $[\Phi_2, \Psi_2] \subset \mathbb{R}$ such that $a \in [\Psi_1, \Phi_1]$ iff the left endpoint of $Rot(F_a)$ is $\frac{p}{q}$, and $a \in [\Phi_2, \Psi_2]$ iff the right endpoint of $Rot(F_a)$ is $\frac{p}{q}$. Moreover, the following hold:
(i) $\Phi_2 \leq \Phi_1$, $\Phi_2 \leq \Psi_1$, and $\Psi_2 \leq \Phi_1$,
(ii) each of $[\Psi_1, \Phi_1]$ and $[\Phi_2, \Psi_2]$ are nondegenerate closed intervals,
(iii) $a = \Phi_1$ ($\Phi_2$) iff $F_a^q(x) \geq x + p$ ($\leq x + p$) and equality holding for some $x$,
(iv) $a = \Psi_1$ ($\Psi_2$) iff $F_l^q(x) \leq x + p$ ($F_u^q(x) \geq x + p$) and equality holding for some $x$,
(v) if $\Psi_2 < \Phi_1$ and $S(F) < 0$, then $a = \Psi_2$ iff $F_{a;u}^q(x) \geq x + p$ and equality holds at $C^+$,
(vi) if $\Phi_2 < \Psi_1$ and $S(F) < 0$, then $a = \Psi_1$ iff $F_{a;l}^q(x) \leq x + p$ and equality holds at $K^-$,
(vii) $Rot(F_a) = \{\frac{p}{q}\} \iff a \in [\Psi_1, \Psi_2]$.

If $\{F_{a,b}\}$ is a family of degree-one lifts satisfying the conditions for Class $\mathcal{N}$ except possibly the negative Schwarzian condition, then $\Psi_1, \Psi_2, \Phi_1$, and $\Phi_2$ are continuous functions of $b$. See Figure 2(b). For more discussion of these functions see [B, MT].

We now describe the Farey organization of the rational numbers. (See, e.g., [GT, HW, LT].) Throughout this note all fractions are in lowest terms and in $[0,1]$ except where otherwise noted. We denote $\mathbb{Q} \cap [0,1]$ by $\mathbb{Q}_{[0,1]}$. The fractions $\frac{p_1}{q_1}$ and $\frac{p_2}{q_2}$ are said to be *Farey neighbors* if they satisfy $p_1 q_2 - p_2 q_1 = \pm 1$. A pair of fractions $\frac{p_1}{q_1}$ and $\frac{p_2}{q_2}$ can be combined to form their *mediant*, defined as $\frac{p_1 + p_2}{q_1 + q_2}$. This is not necessarily in lowest terms, but if $\frac{p_1}{q_1}$ and $\frac{p_2}{q_2}$ are Farey neighbors, $\frac{p_1}{q_1} \oplus \frac{p_2}{q_2} \stackrel{def}{=} \frac{p_1 + p_2}{q_1 + q_2}$ is in lowest terms, and is called their *Farey sum*.

**Remark 1** [GT, HW, LT]. Every rational number $\frac{p}{q} \in (0,1)$ can be expressed uniquely as a Farey sum $\frac{p}{q} = \frac{p_1}{q_1} \oplus \frac{p_2}{q_2}$.

The fraction $\frac{p_1}{q_1}$ is called the *left parent* of $\frac{p}{q}$, and $\frac{p_2}{q_2}$ the *right parent*. We call the fractions in $(0,1)$ whose left parent is $\frac{p}{q}$ the *right children* of $\frac{p}{q}$, and those whose right parent is $\frac{p}{q}$ the *left children* of $\frac{p}{q}$. (The fraction $\frac{0}{1}$ has only right children, and $\frac{1}{1}$ only left children.)

**Definition 1.** Let $\frac{p}{q} = \frac{p_1}{q_1} \oplus \frac{p_2}{q_2}$ be given with $\frac{p}{q} \in (0,1)$. For each $j \geq 0$ set

$$l_j(\frac{p}{q}) = \frac{jp + p_1}{jq + q_1} \text{ and } r_j(\frac{p}{q}) = \frac{jp + p_2}{jq + q_2}.$$

For $\frac{0}{1}$ we define only the sequence $\{r_j\}$ where $\frac{p_2}{q_2} = \frac{1}{1}$, and for $\frac{1}{1}$ we define only the sequence $\{l_j\}$ where $\frac{p_1}{q_1} = \frac{1}{1}$.

The $l_j(\frac{p}{q})$'s, $j = 1, 2, ...$, are the left children of $\frac{p}{q}$, and the $r_j(\frac{p}{q})$'s, $j = 1, 2, ...$, are its right children. Observe that the $l_j$'s converge up to $\frac{p}{q}$ and that the $r_j$'s converge down to $\frac{p}{q}$. Note also that $r_0$ is the right parent, and $l_0$ is the left parent.

The *Farey tree* is the infinite tree whose vertices are the set of rational numbers in $[0,1]$ with one edge between every fraction $\frac{p}{q}$ and each of its first children ($j = 1$). See Figure 4(v). The vertices fall into levels indexed by the non-negative integers. Vertices $\frac{0}{1}$ and $\frac{1}{1}$ are at level 0, while the vertices at level $n > 0$ are the children of the vertices at level $n - 1$. For example, $\frac{1}{3}$ is at level 2, and $\frac{4}{7}$ is at level 4. If two fractions, $\frac{p_m}{q_m}$ at level $m$ and $\frac{p_n}{q_n}$ at level n, are joined by a path in the tree that is monotone with respect to the levels, we say $\frac{p_m}{q_m}$ is *higher* than $\frac{p_n}{q_n}$ if $m < n$. For example, $\frac{1}{3}$ is higher than $\frac{3}{8}$; but $\frac{1}{3}$ is not higher than $\frac{4}{7}$. (as there is no monotone path between $\frac{1}{3}$ and $\frac{4}{7}$). For a generalization of the Farey tree to an *extended Farey tree* and for historical notes on the subject, see [LT].

**Remark 2.** For each irrational $\omega \in (0,1)$ there is precisely one infinite path in the Farey tree for which the values along the path converge to $\omega$. Similarly for $\frac{0}{1}$ and $\frac{1}{1}$. Whereas for each rational $\frac{p}{q} \in (0,1)$ there are exactly two infinite paths in the Farey tree which converge to $\frac{p}{q}$. These paths, which correspond to the continued fraction truncations, are $\{l_j\}_{j \geq 0}$ and $\{r_j\}_{j \geq 0}$ (as defined in Definition 1).

We now define the *Farey web* for a family of lifts $\{F_{a,b}\}$ in Class $\mathcal{N}$.

**Definition 2.** Let $\{F_{a,b}\}$ be a family of lifts in Class $\mathcal{N}$. For each $\frac{p}{q} \in \mathbb{Q}_{[0,1]}$, define $\mathcal{R}_{\frac{p}{q}} = \{(a,b) \mid F_{a,b}^q(K^-) = C^+ + p, F_{a,b}^i(K^-) \notin (K^-, K] + \mathbb{Z} \text{ for } 0 < i \leq q, F_{a,b}^i(K^-) \neq C^+ + p \text{ for } 0 < i < q\}$ and define $\mathcal{L}_{\frac{p}{q}} = \{(a,b) \mid F_{a,b}^q(C^+) = K^- + p, F_{a,b}^i(C^+) \notin [C, C^+) + \mathbb{Z} \text{ for } 0 < i \leq q, F_{a,b}^i(C^+) \neq K^- + p \text{ for } 0 < i < q\}$. We relax the second and third restrictions in each case if $K^- = K = C = C^+$, i.e., on the critical line. (Again, note that to avoid clutter we have suppressed subscripts that would indicate the dependence of $K^-$ and $C^+$ on the parameter $b$.) The *Farey web* for the family $\{F_{a,b}\}$ is the collection of loci given by $\mathcal{R}_{\frac{p}{q}}$ and $\mathcal{L}_{\frac{p}{q}}$ as $\frac{p}{q}$ varies over $\mathbb{Q}_{[0,1]}$. (For rotation numbers in $[0,1]$ $\mathcal{L}_{\frac{0}{1}}$ and $\mathcal{R}_{\frac{1}{1}}$ do not concern us.)

**Theorem 1.** *Let $\{F_{a,b}\}$ be a family of lifts satisfying the conditions for Class $\mathcal{N}$ except possibly the negative Schwarzian condition, and let $\frac{p}{q} \in \mathbb{Q}_{[0,1]}$ be given. Then $\mathcal{R}_{\frac{p}{q}}$ and $\mathcal{L}_{\frac{p}{q}}$ each intersect every horizontal line in parameter space exactly once. Moreover, if for each $b$ we let $\mathcal{R}_{\frac{p}{q}}(b)$ and $\mathcal{L}_{\frac{p}{q}}(b)$ denote the intersection of $\mathcal{R}_{\frac{p}{q}}$ and $\mathcal{L}_{\frac{p}{q}}$ with the horizontal line at $b$ in parameter space, then $\mathcal{R}_{\frac{p}{q}}(b)$ and $\mathcal{L}_{\frac{p}{q}}(b)$ are continuous functions of $b$. $\mathcal{L}_{\frac{p}{q}}$ and $\mathcal{R}_{\frac{p}{q}}$ both intersect the critical line in the unique point, $B_{\frac{p}{q}}$, of the critical line where the critical point is $\frac{p}{q}$-periodic.*

**Theorem 2.** *Let $\frac{p}{q} \in \mathbb{Q}_{[0,1]}$ be given with $\frac{p}{q} = \frac{p_1}{q_1} \oplus \frac{p_2}{q_2}$, and let the sequences $\{l_j\}_{j \geq 0}$ and $\{r_j\}_{j \geq 0}$ be as Definition 1. Let $\{F_{a,b}\}$ be a family of lifts in Class $\mathcal{N}$. Let $T_{\frac{p}{q}}$ denote a tip of $Lock_{\frac{p}{q}}$, and for each $j \geq 0$, let $T_{l_j}$ and $T_{r_j}$ denote tips of $Lock_{l_j}$ and $Lock_{r_j}$ respectively. Then for each $j \geq 0$, $T_{l_j} \in \mathcal{L}_{\frac{p}{q}}$, $T_{r_j} \in \mathcal{R}_{\frac{p}{q}}$, $T_{\frac{p}{q}} \in \mathcal{R}_{l_j}$, and $T_{\frac{p}{q}} \in \mathcal{L}_{r_j}$.*

Figure 3, a sketch of some of the strands of the web associated with $Lock_{\frac{p}{q}}$, illustrates this Theorem. The Farey web serves to link together all the special points of the parameter-plane ($B_{\frac{p}{q}}$'s and $T_{\frac{p}{q}}$'s) pertaining to cycles which are "twist", i.e., combinatorially the same as a cycle of a rigid rotation. (See Definition 3 below.) cycles. A web with a similar character in the context of bimodal interval maps has been described in [RS1, RS2, RS3, RT].

**Theorem 3.** *Let $\{F_{a,b}\}$ be a family in Class $\mathcal{N}$. Fix $b_0$ in parameter space, fix $\frac{p}{q} \in \mathbb{Q}_{[0,1]}$ with $\frac{p}{q} = \frac{p_1}{q_1} \oplus \frac{p_2}{q_2}$, and let $\{r_j\}$ and $\{l_j\}$ be as in Definition 1. For each $j \geq 0$ let $L_j$ denote the first coordinate of the point in $\mathcal{L}_{r_j}$ which lies on the horizontal line $b = b_0$; similarly let* $R_j$ *denote the first coordinate of the point in* $\mathcal{R}_{l_j}$ *which lies on the horizontal line* $b = b_0$. *Set* $\Psi_1 = \Psi_1(b_0)$ *and* $\Psi_2 = \Psi_2(b_0)$. *Then, exactly one of the following hold.*

(1) $\quad L_0 > L_1 > \cdots > L_j > \cdots > \Psi_2 > \Psi_1 > \cdots > R_j > \cdots > R_1 > R_0$

(2) $\quad L_0 < L_1 < \cdots < L_j < \cdots < \Psi_2 < \Psi_1 < \cdots < R_j < \cdots < R_1 < R_0$

(3) $\quad L_0 = L_1 = \cdots = L_j = \cdots = \Psi_1 = \Psi_2 = \cdots = R_j = \cdots = R_1 = R_0$

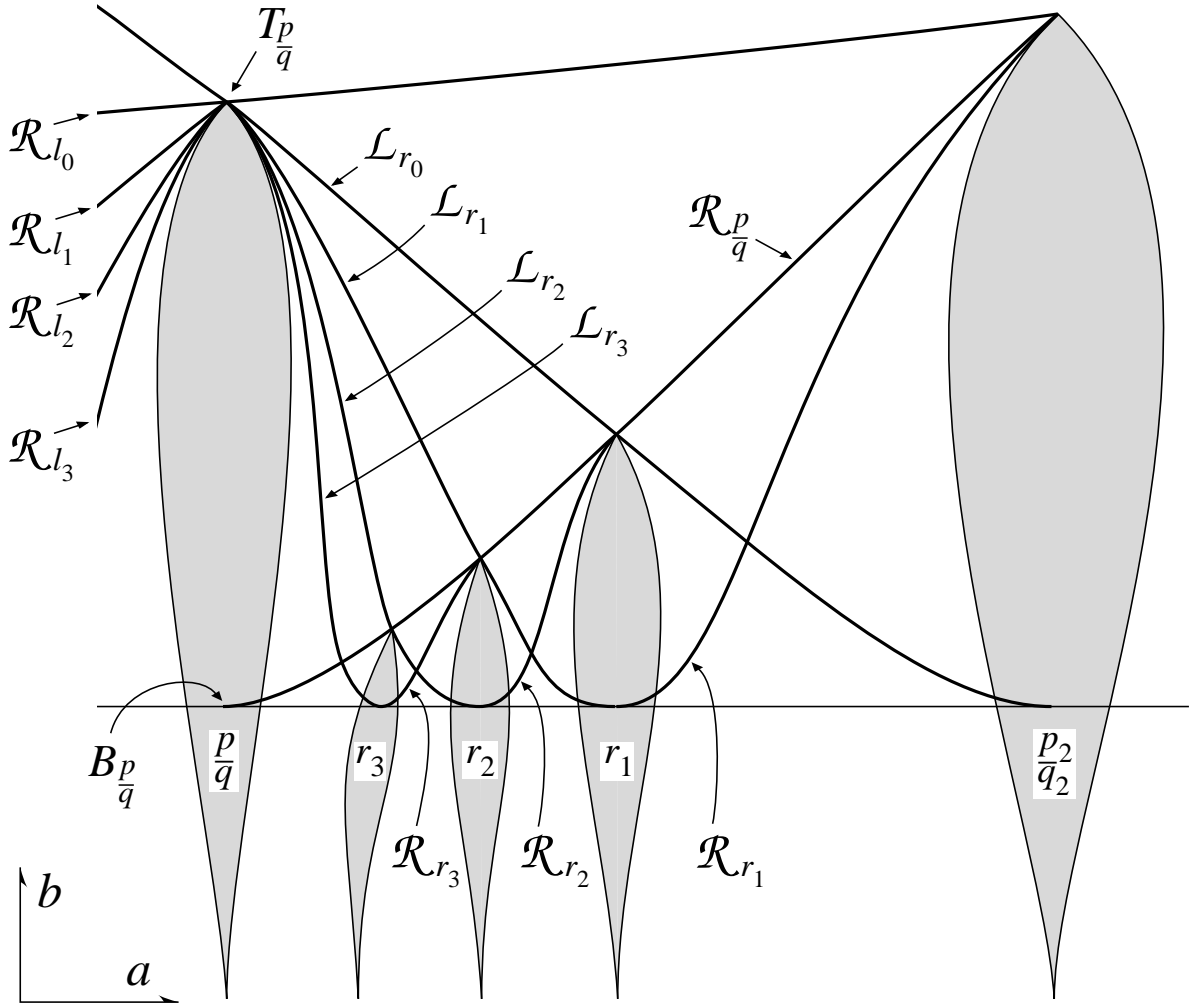

**Figure 3.** Locking regions for $\frac{p}{q}$, its right parent $\frac{p_2}{q_2}$, and three of its right children, linked together by strands of the Farey web. (All $r_j$'s and $l_j$'s are those of $\frac{p}{q}$.)

A consequence of Theorems 1, 2, and 3 is that the simplest Farey web for a family in class $\mathcal{N}$ is a homeomorphic image of the infinite graph whose construction we depict in Figure 4 and describe next. It is a plane graph in a rectangle. Construct a (middle-third, say) Cantor set on the base of the rectangle, leaving room for an extra "level-0" hole on either side. The graph has a pair of vertices, $T_{\frac{p}{q}}$ and $B_{\frac{p}{q}}$, for each rational $\frac{p}{q} \in [0,1]$. Every $B_{\frac{p}{q}}$ is on the base: if $\frac{p}{q}$ appears at the $n^{\text{th}}$ level of the Farey tree, then the common first coordinate of $B_{\frac{p}{q}}$ and $T_{\frac{p}{q}}$ is the center of one of the intervals that is removed at the $n^{\text{th}}$ stage of the hole-punching construction of the Cantor set.

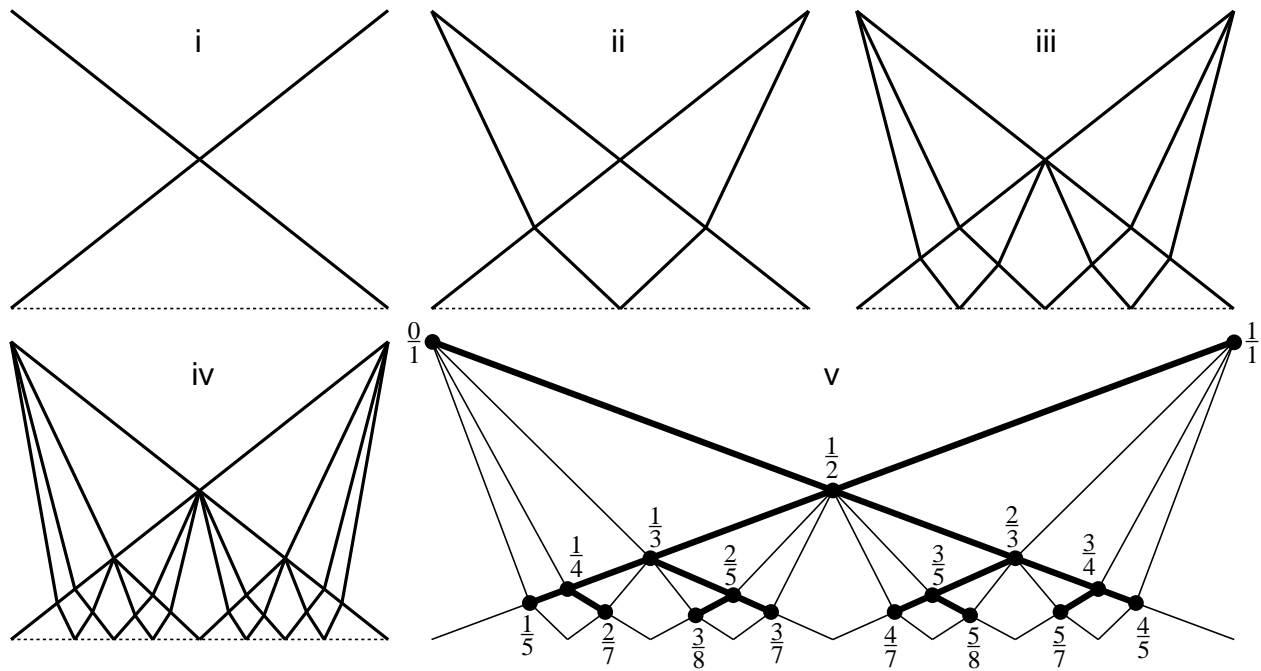

Figure 4. Construction of the Farey web, and its relation to the Farey tree.

Start by placing the level-0 vertices, $T_{\frac{0}{1}}$ and $T_{\frac{1}{1}}$, on the top edge of the rectangle, and draw the segment $T_{\frac{0}{1}} B_{\frac{1}{1}}$. The second coordinate of the remaining $T_{\frac{p}{q}}$'s is determined implicitly by the construction.

At each stage of the construction, take in turn each pair of vertices on the base that are adjacent upon completion of the previous stage. Call them $B_{\frac{p_1}{q_1}}$ and $B_{\frac{p_2}{q_2}}$, choosing the indices (1 and 2) such that $T_{\frac{p_1}{q_1}}$ and $B_{\frac{p_2}{q_2}}$ are already linked by an edge (not necessarily $\frac{p_1}{q_1} < \frac{p_2}{q_2}$). The position of $T_{\frac{p}{q}}$, $\frac{p}{q} = \frac{p_1}{q_1} \oplus \frac{p_2}{q_2}$, is determined by its specified first coordinate as described above and the requirement that it lie on the segment $T_{\frac{p_1}{q_1}} B_{\frac{p_2}{q_2}}$. Draw the "dogleg" from $T_{\frac{p_2}{q_2}}$ through $T_{\frac{p}{q}}$ to $B_{\frac{p_1}{q_1}}$. Doing this for every pair completes one stage of the construction.

The relationship with the Farey tree is indicated in Figure 4(v). In the web, each $T_{\frac{p}{q}}$ is joined by a strand to both of its "parents", $T_{\frac{p_1}{q_1}}$ and $T_{\frac{p_2}{q_2}}$, instead of only one of them as in the tree. In the correspondence with the Farey web in the parameter plane, $T_{\frac{p}{q}}$ corresponds the tip of the locking region $Lock_{\frac{p}{q}}$, and $B_{\frac{p}{q}}$ to the unique point on the critical line where the critical point belongs to a $\frac{p}{q}$-cycle. Each dogleg corresponds to a curve $\mathcal{L}$ or $\mathcal{R}$ with the portion above $T_{\frac{p_2}{q_2}}$ omitted. Figure 1 shows the result of carrying out eight stages of the construction.

Combining our Theorems 1, 2, and 3 with the result of [EKT], we can conclude that in the standard sine family the Farey web is a homeomorphic image of the result of the above construction. However, in class $\mathcal{N}$, multiple coincidences of a compatible $\mathcal{L}$ and $\mathcal{R}$ can occur, as in the sketch, Figure 5, and this can give rise to multiple components of a $Lock_{\frac{p}{q}}$ as shown.

Nevertheless, the regularity of the web characterized by Theorems 1 and 3, does impose significant restrictions on the organization of the frequency-locking regions. Examples are the following theorems and corollaries. Both Theorem 4 (in slightly weaker form) and Corollary 1 are proven in [BT]. We use the Farey web (in particular Theorem 3) to give a short proof for Theorem 4.

**Theorem 4.** *Let $\{F_{a,b}\}$ be a family of lifts in Class $\mathcal{N}$. Fix $b_0$, above the critical line, and consider the one-parameter family $\{F_{a,b_0}\}_{a \in \mathbb{R}}$. Let $\frac{p}{q}$ be higher than $\frac{p'}{q'}$ on the Farey tree. If there exists an $a$ such that $F_{a,b_0}$ is $\frac{p'}{q'}$-locked, then the line $b = b_0$ intersects the $\frac{p}{q}$-locking region in a segment of positive length.*

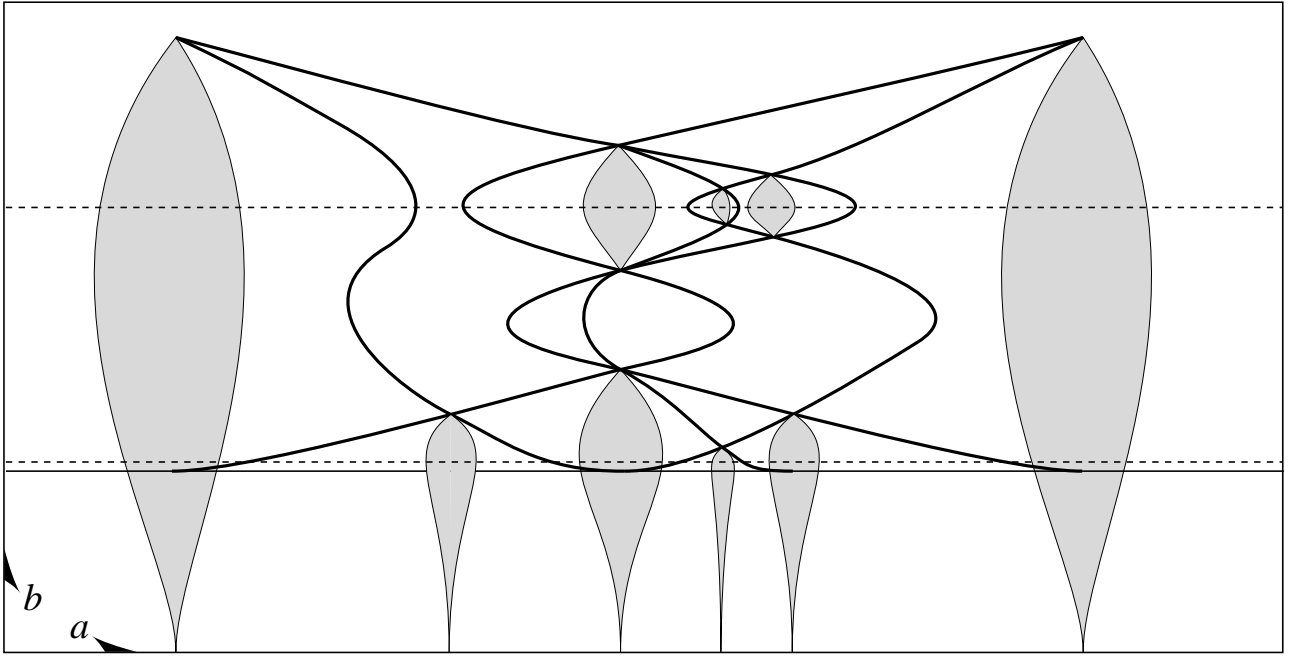

**Figure 5.** On a horizontal line, $\frac{p'}{q'}$-locking implies locking for all $\frac{p}{q}$ higher than $\frac{p'}{q'}$ in the Farey tree.

**Corollary 1.** *Let $\{F_{a,b}\}$ be a family of lifts in Class $\mathcal{N}$. If $\frac{p}{q}$ is higher than $\frac{p'}{q'}$ in the Farey tree, then the reduced height of $Lock_{\frac{p}{q}}$ is greater than the reduced height of $Lock_{\frac{p'}{q'}}$.*

Of course, Corollary 1 is immediate from Theorem 4, but can also be proven, independent of Theorem 4, with Theorem 3. Theorem 5 (and thereby Corollary 2) are proven in [FT] for $\frac{p}{q} = \frac{0}{1}$ only. We generalize the arguments in [FT] for the $\frac{0}{1}$ tongue to an arbitrary one. We use Corollary 1 in our proof of Theorem 5.

**Theorem 5.** *Let $\frac{p}{q}$ be given and let $\{F_{a,b}\}$ be a family of lifts in Class $\mathcal{N}$. For each $j \geq 1$, let $r_j$ and $l_j$ be as in Theorem 2. Then there exists $\varepsilon_{\frac{p}{q}} > 0$ such that all the tips $T_{r_j}$ and $T_{l_j}$* are bounded away from the critical line by a distance at least $\varepsilon_{\frac{p}{q}}$. Furthermore, the sequences of tips $\{T^I_{l_j}\}$ and $\{T^I_{r_j}\}$ each converge, monotonically in b, to a point on the boundary of the $\frac{p}{q}$ locking region above the critical line. On the other hand, for any $\omega \in \mathbb{R} \setminus \mathbb{Q}$, let $\{\frac{c_j}{d_j}\}_{j \geq 1}$ be the path in the Farey tree which converges to $\omega$. Then the $T^I_{\frac{c_j}{d_j}}$'s converge to the critical line with $j$.

(Observe how in Figure 3, the tips for the right children of $\frac{p}{q}$ all lie on $\mathcal{R}_{\frac{p}{q}}$ which intersects the boundary of $Lock_{\frac{p}{q}}$ above the critical line.)

**Corollary 2.** *For families $\{f_{a,b}\}$ in Class $\mathcal{N}$, the Arnold tongue for every rational number has on each side an arc of its boundary at whose points the boundary of zero-width rotation interval for the family is not locally connected.*

## 3. Preparatory material on degree-one lifts

For the convenience of the reader and for completeness we collect some basic facts about degree-one lifts and some known facts about degree-one lifts in Class $\mathcal{B}$. We use these facts throughout Section 4.

**Fact 2** [ALM, Lemma 3.7.12]**.** The function $F \mapsto \rho(F)$ is continuous on the space of all non-decreasing degree-one lifts.

**Fact 3** [ALM, proof Lemma 3.7.19]**.** Let $F$ and $G$ be non-decreasing degree-one lifts. If $F < G$, then $F^n < G^n$ for all $n \in \mathbb{N}$.

**Fact 4** [ALM, Theorem 4.7.20]**.** Let $F$ be a degree-one lift. Then $Rot(F) = [\rho_{F_l}, \rho_{F_u}]$.

**Remark 3.** Let $F$ be a degree one lift and consider the one-parameter family $\{F_a = F + a\}_{a \in \mathbb{R}}$. For ease of notation, set $F_{a;l} = (F_a)_l$ and $F_{a;u} = (F_a)_u$. For every $a \in \mathbb{R}$ we have that $F_{a;l} = F_l + a$ and $F_{a;u} = F_u + a$. Hence, it follows from Fact 3 that $a_1 < a_2$ implies both $\rho(F_{a_1;l}) \leq \rho(F_{a_2;l})$ and $\rho(F_{a_1;u}) \leq \rho(F_{a_2;u})$.

**Definition 3.** Let F be a degree-one lift and let $x \in \mathbb{R}$. The *lift-orbit of $x$* is defined to be $orb(x, F) + \mathbb{Z}$, where $orb(x, F) = \{x, F(x), F^2(x), \cdots\}$. We say that the lift-orbit of x is a *twist lifted orbit* provided that $F$ restricted to $orb(x, F) + \mathbb{Z}$ is non-decreasing. We say that $orb(x, F) + \mathbb{Z}$ is a *lifted cycle* provided there is some $n \in \mathbb{N}$ such that $F^n(x) - x \in \mathbb{Z}$. The lift-orbit $orb(x, F) + \mathbb{Z}$ is called a *lifted m-cycle* if $orb(x, F) + \mathbb{Z}$ is a lifted cycle and if $m = min\ \{n \in \mathbb{N} \mid F^n(x) - x \in \mathbb{Z}\}$. Lastly, $orb(x, F) + \mathbb{Z}$ is a *twist lifted m-cycle* if $orb(x, F) + \mathbb{Z}$ is a lifted $m$-cycle and a twist lifted orbit (in this case we say that $x$ or $F$ has a twist lifted $m$-cycle) [ALM, sections 3.2, 3.3 and 3.7].

Fact 5 is well known and for convenience we cite [BT]. This fact describes the combinatorics of a twist lifted $\frac{p}{q}$-cycle (*i.e.* a twist lifted $q$-cycle with $F^q(x) = x + p$).

**Fact 5** [BT, Lemma 4]**.** Let $F$ be a degree-one lift and let $x \in \mathbb{R}$ have a twist lifted $q$-cycle with $F^q(x) - x = p$. Set, $orb(x, F) + \mathbb{Z} = \{\cdots, y_{-2}, y_{-1}, y_0, y_1, y_2, \cdots\}$ with $\cdots < y_{-2} < y_{-1} < y_0 < y_1 < y_2 < \cdots$ and with $(orb(x, F) + \mathbb{Z}) \cap [0, 1) = \{y_0, y_1, \cdots y_{q-1}\}$. Recall that $\frac{p}{q} = \frac{p_1}{q_1} \oplus \frac{p_2}{q_2}$. Lastly, let $i_0 \in \{0, 1, \cdots, q-2\}$ and $i_1 \in \{1, 2, \cdots q-1\}$. Then, $F^{q_1}(y_{i_0})\ mod\ 1 = y_{i_0+1}$, $F^{q_1}(y_{i_0}) \in [p_1, p_1+1)$, $F^{q_2}(y_{i_1})\ mod\ 1 = y_{i_1-1}$, and $F^{q_2}(y_{i_1}) \in [p_2, p_2+1)$.

**Fact 6** [MT, Proposition 4.2]. Let $F$ be a degree-one lift which has precisely one critical point and negative Schwarzian elsewhere, or precisely two critical points, both turning points, and negative Schwarzian elsewhere. Let $\frac{p}{q} \in Rot(F)$. Then either $F$ has either precisely two twist lifted $\frac{p}{q}$-cycles and the graph of $F^q(x) - p$ crosses the diagonal transversally and in opposite directions at their points, or $F$ has precisely one twist lifted $\frac{p}{q}$-cycle and the graph of $F^q(x) - p$ does not cross the diagonal at its points. Moreover, the same statements hold for any $F_l$ with $Rot(F_l) = \frac{p}{q}$ or any $F_u$ with $Rot(F_u) = \frac{p}{q}$.

**Fact 7** [BT, Lemma 6]. Let $F$ be a degree-one lift in Class $\mathcal{B}$ with precisely two turning points. Let $\frac{p}{q} \in \mathbb{Q}_{[0,1]}$ and suppose that $Rot(F) = \{\frac{p}{q}\}$. Then, F has at least one twist lifted $\frac{p}{q}$-cycle with two consecutive points $M < N$ such that $M \in [0, K^-]$ and $N \in [C^+, 1)$. As is pointed out in [BT] this fact is known, for convenience we simply cite [BT].

**Fact 8** [MT, Lemma 4.3]. Let $F$ be a degree-one lift such that $F$ is either increasing or has at most one interval $[C, K] \subset [0, 1]$ where F is non-increasing. If $F$ has two twist lifted $\frac{p}{q}$-cycles, then each has precisely one point between consecutive points of the other.

**Fact 9** [B, Lemma 1.4]. Let $\{F_{a,b}\}$ be a family of degree one lifts and let $\frac{p}{q} \in \mathbb{Q}$. Then, the bifurcation that occurs when $\frac{p}{q}$ is lost from $Rot(F_{a,b})$ is a $\frac{p}{q}$-tangent bifurcation.

**Fact 10** [B, proof of Theorem 4.1]. Let $\{F_{a,b}\}$ be a family of degree one lifts satisfying the conditions for Class $\mathcal{N}$ except possibly the negative Schwarzian condition. Let $\frac{p}{q} \in \mathbb{Q}$. Then, there is a unique point on the critical line and within the $\frac{p}{q}$ Arnold tongue (*i.e.*, not on a boundary) such that $K_b = C_b$ is part of a twist $\frac{p}{q}$-cycle for $F_{a,b}$.

Fact 11 follows from work in [B,MT] and is known; however we could not find an instance where Fact 11 is explicitly stated, so for completeness we provide a proof of it in Section 4.

**Fact 11.** Let $\{F_{a,b}\}$ be a family of lifts in Class $\mathcal{N}$. Let $\frac{p}{q}$ be given. Then for each $b$ we have the following:
i) if $\Psi_1(b) = \Psi_2(b)$, then for $F = F_{\Psi_1(b),b}$ we have that $K_F^-$ and $C_F^+$ are consecutive points in a twist lifted $\frac{p}{q}$-cycle for $F$, and conversely
ii) if $K_F^-$ and $C_F^+$ are consecutive points in a twist lifted $\frac{p}{q}$-cycle for some $F = F_{a,b}$, then $a = \Psi_1(b) = \Psi_2(b)$.

**Remark 4.** Let $\{F_{a,b}\}$ be a family in Class $\mathcal{N}$, let $\frac{p}{q} = \frac{p_1}{q_1} \oplus \frac{p_2}{q_2}$ with $\frac{p}{q} \in \mathbb{Q}_{[0,1]}$, and let $F$ be the lift corresponding to a tip of $Lock_{\frac{p}{q}}$. From Facts 4 and 10 we have that $K^-$ and $C^+$ are consecutive points in a twist lifted $\frac{p}{q}$-cycle for $F$. Hence, $F^{q_1}(K^-) = C^+ + p_1$, $F^{q_2}(C^+) = K^- + p_2$, and this twist lifted cycle misses $(K^-, C^+) + \mathbb{Z}$.

Figure 6, which shows iterates of a lift of the sine map at the tip of $Lock_{\frac{3}{8}}^I$, illustrates this remark.

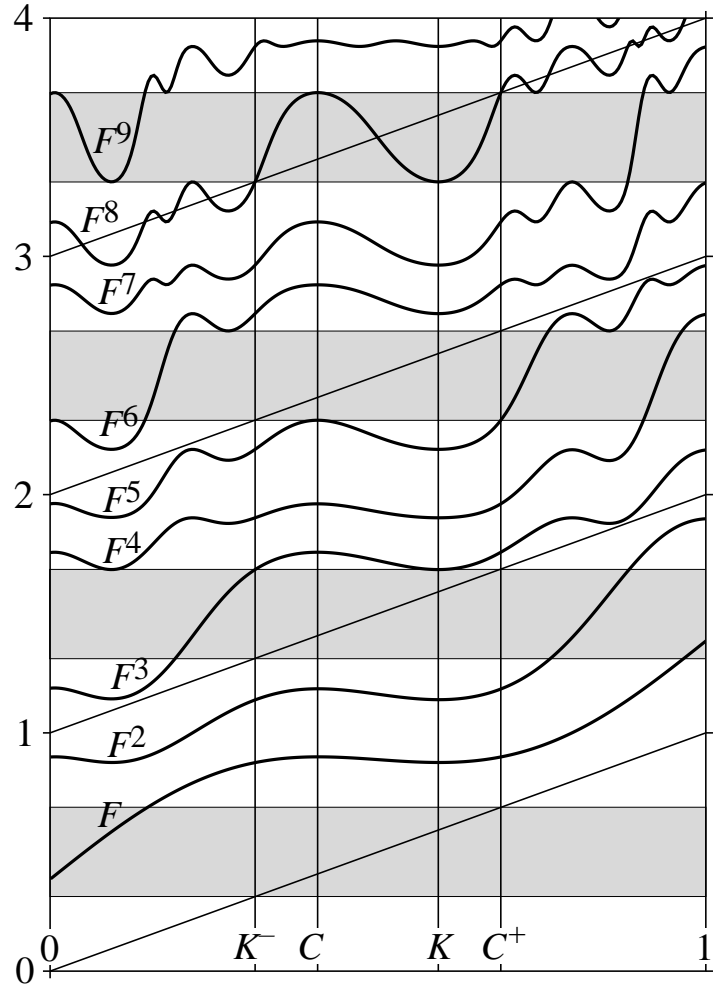

**Figure 6.** At $T^I_{\frac{3}{8}}$ in the standard sine family.

## 4. Proofs

*Proof of Fact 11.* (i) Assume that $\Psi_1(b) = \Psi_2(b)$ and set $F = F_{\Psi_1(b),b}$. We will show that $K^-$ and $C^+$ are consecutive points in a twist lifted $\frac{p}{q}$-cycle. Fact 1 (i)-(ii) and $\Psi_1(b) = \Psi_2(b)$ imply that $\Phi_2(b) < \Psi_1(b) = \Psi_2(b) < \Phi_1(b)$. Hence, parts (v) and (vi) of Fact 1 hold for $F$ and therefore $C^+$ and $K^-$ are each part of a twist lifted $\frac{p}{q}$-cycle for $F$ missing $[C, C^+) + \mathbb{Z}$ and $(K^-, K] + \mathbb{Z}$ respectively. Hence if $C^+$ and $K^-$ are part of same twist lifted $\frac{p}{q}$-cycle then they are consecutive points in the cycle. Thus, if $F$ has precisely one twist lifted $\frac{p}{q}$-cycle we are done. From Fact 6, $F$ has at most two twist lifted $\frac{p}{q}$-cycles. If $F$ has two twist lifted $\frac{p}{q}$-cycles, then, by Fact 6, the graph of $F^q$ crosses '$x + p$' at both $K^-$ and $C^+$ and therefore, use Fact 1 (v)-(vi), the graph of $F^q$ crosses '$x + p$' in the same direction at both $K^-$ and at $C^+$. Thus, by Fact 6, $K^-$ and $C^+$ are in the same twist lifted $\frac{p}{q}$-cycle and we are done.

(ii) Set $F = F_{a,b}$ and assume that $K^-$ and $C^+$ are consecutive points in a twist lifted

$\frac{p}{q}$-cycle (say $\mathcal{O}$) for $F$. We show that $a = \Psi_1(b) = \Psi_2(b)$. First notice that $\mathcal{O}$ is a twist lifted $\frac{p}{q}$-cycle for $F_l$, $F_u$, and $F$ since $K^-$ and $C^+$ are consecutive points in $\mathcal{O}$. Thus, by Fact 1 (iv), it suffices to show that $F_l^q(x) \leq x + p$ and $F_u^q(x) \geq x + p$ for all $x$. We show that $F_l^q(x) \leq x + p$ for all $x$, as the other inequality is similar.

Suppose to the contrary that $F_l^q(x) \leq x + p$ does not hold for all $x$. Then, by Fact 6, $F_l$ has two twist lifted $\frac{p}{q}$-cycles, say $\mathcal{O}$ (as before containing $K^-$ and $C^+$) and $\mathcal{O}_1$; hence, again using Fact 6, $F_l^q(x)$ crosses '$x + p$' at $K^-$ from above to below and therefore $F^q(x)$ also crosses '$x + p$' at $K^-$ from above to below (use $F^q(x) = F_l^q(x)$ for $x \leq K^-$ and sufficiently close to $K^-$). Since $K^-$ and $C^+$ are consecutive points in $\mathcal{O}$ for $F$, we have that $F^q$ also crosses '$x + p$' at $C^+$ from above to below; thus $F_u^q(x)$ crosses $x + p$ at $C^+$ from above to below (again, use $F^q(x) = F_u^q(x)$ for $x \geq C^+$ and sufficiently close to $C^+$). Therefore $F_u$ (by Fact 6) also has two twist lifted $\frac{p}{q}$ cycles; say $\mathcal{O}$ (as before containing $K^-$ and $C^+$) and $\mathcal{O}_2$. Notice that all three cycles, $\mathcal{O}$, $\mathcal{O}_1$, and $\mathcal{O}_2$ are twist lifted $\frac{p}{q}$-cyles for $F$ and therefore, by Fact 6, cannot all be distinct. By construction, $\mathcal{O} \neq \mathcal{O}_1$ and $\mathcal{O} \neq \mathcal{O}_2$; hence $\mathcal{O}_1 = \mathcal{O}_2$. However, $\mathcal{O}_1$ misses $[K^-, K]$ and $\mathcal{O}_2$ misses $[C, C^+]$, since they are cycles for $F_l$ and $F_u$ respectively and not $\mathcal{O}$. Fact 8 (applied to $\mathcal{O}$ and $\mathcal{O}_1$) implies that $\mathcal{O}_1$ intersects $(K^-, C^+)$ and therefore (use $\mathcal{O}_1$ misses $[K^-, K]$) $\mathcal{O}_1$ intersects $(K, C^+)$. Hence $\mathcal{O}_1 \neq \mathcal{O}_2$, since $\mathcal{O}_2$ misses $[C, C^+]$, contradicting $\mathcal{O}_1 = \mathcal{O}_2$. □

*Proof of Theorem 1.* With the notation of Theorem 1, we treat the case of $\mathcal{R}_{\frac{p}{q}}$; the case of $\mathcal{L}_{\frac{p}{q}}$ is similar. We want to show that every horizontal line in parameter space (on or above the critical line) intersects $\mathcal{R}_{\frac{p}{q}}$ exactly once. Notice that if $(a, b) \in \mathcal{R}_{\frac{p}{q}}$, then $F_{(a,b),l}^q(K^-) = C^+ + p$. Fix the parameter $b$ and consider the family of nondecreasing lifts $\{F_{(a,b),l} \mid a \in \mathbb{R}\}$. Fact 3 and the parameter $a$ being a translation parameter imply that there is a unique $a_\star \in \mathbb{R}$ such that $F_{(a_\star,b),l}^q(K^-) = C^+ + p$. Moreover, $F_{(a_\star,b),l}^i(K^-) \notin (K^-, K] + \mathbb{Z}$ for $0 \leq i < q$. For otherwise $F_{(a_\star,b),l}(K^-)$ is periodic under $F_{(a_\star,b),l}$, which is a contradiction; since $C^+ + p$ in this cycle implies that the rotation number for this cycle is of the form $\frac{0}{m}$ with $m > 1$ (a contradiction), or $C^+ + p$ not in the cycle implies that $F_{(a_\star,b),l}(K^-)$ cannot get to $C^+ + p$ as desired.

The second part of the theorem follows from the implicit function theorem and the facts that $K^-$ and $C^+$ are independent of $a$ and $F_{(a,b),l}^q(x)$ is a strictly increasing continuous function of $a$. The last sentence of the theorem follows from Fact 10. □

*Proof of Theorem 2.* Theorem 2 follows from Facts 1, 5, and 11. See also Remark 4. □

*Proof of Theorem 3.* Assume the hypotheses and definitions of Theorem 3. If $\Psi_1 = \Psi_2$, then Facts 5 and 11 imply that (3) of Theorem 3 holds. Hence, we assume that $\Psi_1 \neq \Psi_2$ and show that exactly one of (1) or (2) of Theorem 3 hold. To do this, we prove the following hold for each $j \geq 0$:

(4) $\qquad\qquad\qquad\mathsf{R}_j < \Psi_1 \iff \Psi_1 < \Psi_2 \iff \Psi_2 < \mathsf{L}_j,$

(5) $\qquad\qquad\qquad\mathsf{R}_j < \mathsf{R}_{j+1} \iff \Psi_1 < \Psi_2 \iff \mathsf{L}_{j+1} < \mathcal{L}_j,$

(6) $\qquad\qquad\qquad\Psi_1 > \Psi_2 \Rightarrow \mathsf{L}_{j+1} \neq \mathsf{L}_j \text{ and } \mathsf{R}_{j+1} \neq \mathsf{R}_j.$

Suppose that (4), (5), and (6) hold. If $\Psi_1 < \Psi_2$, then (4) and (5) imply that (2) of Theorem 3 holds. If $\Psi_1 > \Psi_2$, then (4), (5), and (6) give that $\Psi_2 \geq \mathsf{L}_j$, $\mathsf{L}_{j+1} \neq \mathsf{L}_j$, and $\mathsf{L}_{j+1} \geq \mathsf{L}_j$ for all $j \geq 0$. Hence, if $\Psi_2 = \mathsf{L}_{j+1}$ for some $j$, then $\mathsf{L}_{j+1} = \mathsf{L}_{j+k}$ for all $k \geq 1$, a contradiction.

Thus, if $\Psi_1 > \Psi_2$, then (1) of Theorem 3 holds. Hence, to prove Theorem 3 it suffices to establish (4), (5) and (6). We work with the $\mathsf{L}_j$'s; the arguments involving the $\mathsf{R}_j$'s are similar. Set $F = F_{\Psi_1(b_0),b_0}$, $G = F_{\Psi_2(b_0),b_0}$, and for $j \geq 0$ set $H_j = F_{\mathsf{L}_j,b_0}$. Notice that, from the definition of the $\mathsf{L}_j$'s, for each $j \geq 0$ we have,

$$(7) \qquad H_{j;u}^{jq+q_2}(C^+) = K^- + jp + p_2 \quad \text{and} \quad H_j^{jq+q_2}(C^+) = K^- + jp + p_2.$$

Fix $j \geq 0$.

Assume that $\Psi_1 < \Psi_2$. Hence, use Fact 1(vii), $Rot(G) = \{\frac{p}{q}\}$ and therefore Fact 7 applies to $G$. Let $M$, $N$, be as in Fact 7 for $G$; since $M \in [0, K^-]$ and $N \in [C^+, 1]$, this twist lifted $\frac{p}{q}$-cycle is also a twist lifted $\frac{p}{q}$-cycle for both $G_l$ and $G_u$.

We first show that $\Psi_2 < \mathsf{L}_j$. Suppose to the contrary that $\Psi_1 < \Psi_2$ and $\Psi_2 \geq \mathsf{L}_j$. If $\Psi_2 > \mathsf{L}_j$, then $G_u > H_{j;u}$ and therefore (use Fact 3 for $G_u > H_{j;u}$, (7), $M \leq K^-$, and Fact 5)

$$(8) \qquad G_u^{jq+q_2}(C^+) > H_{j;u}^{jq+q_2}(C^+) = K^- + jp + p_2 \geq M + jp + p_2 = G_u^{jq+p_2}(N).$$

However, (8) implies that $C^+ \neq N$ and thus (8) contradicts $G_u^{jq+p_2}$ being nondecreasing since $N \geq C^+$. If $\mathsf{L}_j = \Psi_2 < \Phi_1$, then Fact 1(v) implies that $C^+$ is part of a twist lifted $\frac{p}{q}$-cycle for $G_u$, say $\hat{\mathcal{O}}$, and thus $G_u = H_{j;u}$ implies that (use (7))

$$(9) \qquad G^{jq+p_2}(C^+) = G_u^{jq+p_2}(C^+) = H_{j;u}^{jq+p_2}(C^+) = K^- + jp + p_2.$$

Hence, (use (9) and Fact 5) $C^+$ and $K^-$ are consecutive points in this twist lifted $\frac{p}{q}$-cycle $\hat{\mathcal{O}}$ and therefore $\Psi_1 = \Psi_2$ by Fact 11. This contradicts $\Psi_1 < \Psi_2$. If $\mathsf{L}_j = \Psi_2 = \Phi_1$, then $M < K^-$. (For otherwise $\Psi_2 = \Phi_1$ implies, using Fact 1(iii), that $G^q(x) \geq x + p$ for all $x$ and hence that $G_l^q(x) \geq x + p$ for all $x$; but $M = K^-$ implies that $G_l^q(K) = G_l^q(K^-) = K^- + p < K + p$.) Thus (use Fact 5, $M < K^-$, (7), and $H_j = G_u$),

$$(10) \qquad G_u^{jq+p_2}(N) = M + jp + p_2 < K^- + jp + p_2 = H_j^{jq+p_2}(C^+) = G_u^{jq+p_2}(C^+).$$

However (10) implies that $N \neq C^+$ and then (10) contradicts $G_u^{jq+p_2}$ being nondecreasing since $N > C^+$. We now have:

$$(11) \qquad \Psi_1 < \Psi_2 \Rightarrow \Psi_2 < \mathsf{L}_j \quad \text{for each} \quad j \geq 0.$$

We now show that $\mathsf{L}_{j+1} < \mathsf{L}_j$ (still assuming that $\Psi_1 < \Psi_2$). Suppose to the contrary that $\mathsf{L}_{j+1} \geq \mathsf{L}_j$. Hence we have (using (11)) $\Psi_1 < \Psi_2 < \mathsf{L}_j \leq \mathsf{L}_{j+1}$. First notice that $\mathsf{L}_j \leq \mathsf{L}_{j+1}$ implies that (use (7))

$$(12) \quad H_{j;u}^{(j+1)q+q_2}(C^+) = H_{j;u}^q(H_{j;u}^{jq+q_2}(C^+)) = H_{j;u}^q(K^- + jp + p_2) = H_{j;u}^q(K^-) + jp + p_2,$$

and therefore (use (12), Fact 3 for $H_{j;u} \leq H_{j+1;u}$, and (7))

$$(13) \quad H_{j;u}^q(K^-) + jp + p_2 = H_{j;u}^{(j+1)q+q_2}(C^+) \leq H_{j+1;u}^{(j+1)q+p_2}(C^+) = K^- + (j+1)p + p_2.$$

From (13) we get,

$$H_{j;u}^{q}(K^{-}) \leq K^{-} + p. \tag{14}$$

However, from Fact 1(iv) we have that,

$$G_{u}^{q}(K^{-}) \geq K^{-} + p. \tag{15}$$

Joining (14) and (15) and applying Fact 3 to $H_{j;u} > G_u$ we get,

$$K^{-} + p \geq H_{j;u}^{q}(K^{-}) > G_{u}^{q}(K^{-}) \geq K^{-} + p,$$

a contradiction. We now have:

$$\Psi_1 < \Psi_2 \Rightarrow \mathsf{L}_{j+1} < \mathsf{L}_j \text{ for each } j \geq 0.$$

We now drop our earlier assumption that $\Psi_1 < \Psi_2$ in order to establish the remaining parts of (4), (5), and (6).

Assume that $\Psi_2 < \mathsf{L}_j$. We show that $\Psi_1 < \Psi_2$. Suppose to the contrary that $\Psi_1 \geq \Psi_2$. If $\Psi_1 = \Psi_2$, then Facts 5 and 11 imply that $\Psi_2 = \mathsf{L}_j$, contradicting our assumption that $\Psi_2 < \mathsf{L}_j$. Hence, suppose $\Psi_1 > \Psi_2$. Then Fact 1 (i)-(ii) imply that $\Phi_2 < \Psi_2 < \Psi_1 < \Phi_1$ and hence (use Fact 1(v)) $C^+$ is part of a twist lifted $\frac{p}{q}$-cycle for $G_u$. Thus,

$$G_{u}^{jq+q_2}(C^+) = G_{u}^{q_2}(C^+) + jp. \tag{16}$$

However $\Psi_2 < \mathsf{L}_j$ gives that $G_u < H_{j;u}$ and therefore that (use Fact 3 and (7)),

$$G_{u}^{jq+q_2}(C^+) < H_{j;u}^{jq+q_2}(C^+) = K^{-} + jp + p_2. \tag{17}$$

Combining (16) and (17) we get

$$G_{u}^{q_2}(C^+) < K^{-} + p_2,$$

and therefore (use Fact 5) the twist lifted $\frac{p}{q}$-cycle of $C^+$ under $G_u$ misses $[K^{-}, C^{+}) + \mathbb{Z}$. Hence this twist $\frac{p}{q}$-cycle is also a twist lifted $\frac{p}{q}$-cycle for $G_l$ and thus $Rot(G) = \{\frac{p}{q}\}$. However, $Rot(G) = \{\frac{p}{q}\}$ implies, (use Fact 1(vii)), that $\Psi_1 \leq \Psi_2$, contradicting our supposition that $\Psi_1 > \Psi_2$. We now have $\Psi_2 < \mathsf{L}_j \Rightarrow \Psi_1 < \Psi_2$ and hence we have:

$$\Psi_1 < \Psi_2 \iff \Psi_2 < \mathsf{L}_j \text{ for each } j \geq 0.$$

Assume that $\mathsf{L}_{j+1} < \mathsf{L}_j$. We show that $\Psi_1 < \Psi_2$. Suppose to the contrary that $\Psi_1 \geq \Psi_2$. If $\Psi_1 = \Psi_2$, then Facts 5 and 11 imply that $\mathsf{L}_{j+1} = \mathsf{L}_j = \Psi_1 = \Psi_2$, contradicting our assumption that $\mathsf{L}_{j+1} < \mathsf{L}_j$. Hence suppose that $\Psi_1 > \Psi_2$. Then, (use Fact 1(vii)) $Rot(H_j) \neq \{\frac{p}{q}\}$ and (use (11)) $\mathsf{L}_j \leq \Psi_2$. Our assumption $\mathsf{L}_{j+1} < \mathsf{L}_j$ gives (use (12), Fact 3 for $H_{j+1;u} < H_{j;u}$, and (7)) that

$$H_{j;u}^{q}(K^{-}) + jp + p_2 = H_{j;u}^{(j+1)q+q_2}(C^+) > H_{j+1;u}^{(j+1)q+q_2}(C^+) = K^{-} + (j+1)p + p_2$$

and hence that

(18) $$H^q_{j;u}(K^-) > K^- + p.$$

We claim the following:

(19) $$H^q_{j;u}(C) < C + p.$$

We first use (19) and then we prove it. Now, (18) and (19) give that the graph of $H^q_{j;u}$ crosses '$x + p$' strictly between $K^-$ and $C$. Thus, by Fact 6, $H_{j;u}$ has two twist lifted $\frac{p}{q}$-cycles, say $\mathcal{O}_1$ and $\mathcal{O}_2$. Neither of these cycles enter $[C, C^+]$ due to (19) and due to $H^q_{j;u}(C^+) = H^q_{j;u}(C)$. Hence both $\mathcal{O}_1$ and $\mathcal{O}_2$ enter $(K^-, C)$; for otherwise at least one cycle would miss $(K^-, C^+) + \mathbb{Z}$ and therefore would also be a twist lifted $\frac{p}{q}$-cycle for $H_{j;l}$ forcing $Rot(H_j) = \{\frac{p}{q}\}$ and contradicting $Rot(H_j) \neq \{\frac{p}{q}\}$. However, both $\mathcal{O}_1$ and $\mathcal{O}_2$ entering $(K^-, C^+)$, (18), (19), and Fact 6 imply that at least one of $\mathcal{O}_1$ or $\mathcal{O}_2$ has two consecutive points in $(K^-, C)$, say $K^- < M < N < C$. Then (use Fact 5, $M > K^-$, and (7)),

$$H^{jq+q_2}_{j;u}(N) = M + jp + p_2 > K^- + jp + p_2 = H^{jq+q_2}_{j;u}(C^+),$$

contradicting $H^{jq+q_2}_{j;u}$ being nondecreasing. Hence, we now have:

$$\mathsf{L}_{j+1} < \mathsf{L}_j \Rightarrow \Psi_1 < \Psi_2 \text{ for each } j \geq 0,$$

modulo (19).

We now prove (19). Suppose to the contrary that $H^q_{j;u}(C) \geq C + p$. Then, using $H^q_{j;u}(C) = H^q_{j;u}(C^+)$, we get that:

(20) $$H^q_{j;u}(C^+) \geq C + p.$$

Since $\Psi_1 > \Psi_2$ we have (use Fact 1(i)-(ii)) that $\Psi_2 < \Phi_1$ and thus (use Fact 1(v)) $G^q_u(C^+) = C^+ + p$. Hence, Fact 3 and $\mathsf{L}_j \leq \Psi_2$ imply that

(21) $$H^q_{j;u}(C^+) \leq G^q_u(C^+) = C^+ + p.$$

Combining (20) and (21) we get that $C + p \leq H^q_{j;u}(C^+) \leq C^+ + p$. Hence, there exists $t \in [C, C^+]$ such that $H^q_{j;u}(t) = t + p$. But, $t \in [C, C^+]$ implies that $H^{jq+q_2}_{j;u}(C^+) = H^{jq+q_2}_{j;u}(t)$ and therefore (use (7))

(22) $$H^{jq+q_2}_{j;u}(t) = H^{jq+q_2}_{j;u}(C^+) = K^- + jp + p_2.$$

Hence, (use Fact 5 and (22)) $K^-$, and $t$ are consecutive points in a twist lifted $\frac{p}{q}$-cycle of $H_{j;u}$ and therefore

(23) $$H^q_{j;u}(K^-) = K^- + p.$$

However, $\mathsf{L}_{j+1} < \mathsf{L}_j$ implies (as for (18)) that $H^q_{j;u}(K^-) > K^- + p$ contradicting (23). This contradiction completes the proof of (19). We now have:

$$\Psi_1 < \Psi_2 \iff \mathsf{L}_{j+1} < \mathsf{L}_j \text{ for each } j \geq 0.$$

It remains to prove (6). Assume that $\Psi_1 > \Psi_2$. As before, (use Fact 1 (i)-(ii)) $\Psi_2 < \Phi_1$ and hence (use Fact 1(v))

$$G_u^q(C^+) = C^+ + p. \tag{24}$$

From (4) and $\Psi_1 > \Psi_2$ we have that

$$\mathsf{L}_j \leq \Psi_2 \text{ for each } j \geq 0. \tag{25}$$

We show that $\mathsf{L}_{j+1} \neq \mathsf{L}_j$. Suppose to the contrary that $\mathsf{L}_{j+1} = \mathsf{L}_j$. If $\mathsf{L}_j = \Psi_2$, i.e., $H_j = G$, then (use (7), $H_{j;u} = G_u$, and (24))

$$K^- + jp + p_2 = H_{j;u}^{q_2}(H_{j;u}^{jq}(C^+)) = G_u^{q_2}(G_u^{jq}(C^+)) = G_u^{q_2}(C^+ + jp) = G_u^{q_2}(C^+) + jp. \tag{26}$$

From (26) we have that $G_u^{q_2}(C^+) = K^- + p_2$ and thus (use (24) and Fact 5) $K^-$ and $C^+$ are consecutive points in a twist lifted $\frac{p}{q}$-cycle for $G_u$ and so also for $G$. Thus, Fact 11 implies $\Psi_1 = \Psi_2$ contradicting our assumption that $\Psi_1 > \Psi_2$. Hence, $\mathsf{L}_j \neq \Psi_2$ and therefore (from (25)) $\mathsf{L}_j < \Psi_2$. By assumption $\mathsf{L}_j = \mathsf{L}_{j+1}$, i.e., $H_j = H_{j+1}$, and thus (use (7))

$$K^- + (j+1)p + p_2 = H_{j;u}^q(H_{j;u}^{jq+q_2}(C^+)) = H_{j;u}^q(K^-) + jp + p_2. \tag{27}$$

From (27) we get,

$$H_{j;u}^q(K^-) = K^- + p;$$

call this twist lifted $\frac{p}{q}$-cycle under $H_{j;u}$, containing $K^-$, $\mathcal{O}$. The definition of $\mathsf{L}_{j+1}$ and $H_{j;u}^{jq+q_2}(C^+) = K^- + jp + p_2$ imply that $\mathcal{O}$ misses $[C, C^+) + \mathbb{Z}$. However, $\Psi_2 < \Psi_1 < \Phi_1$ implies (Fact 1(v)) that $G_u^q(C^+) = C^+ + p$; hence $\mathsf{L}_j < \Psi_2$ gives (use Fact 3 for $H_{j;u} < G_u$)

$$H_{j;u}^q(C^+) < C^+ + p. \tag{28}$$

The cycle $\mathcal{O}$ missing $[C, C^+) + \mathbb{Z}$ and (28) imply that $\mathcal{O}$ misses $[C, C^+] + \mathbb{Z}$. Hence, the orbit of $K^-$ under $H_{j;u}$ must have a point in $(K^-, C) + \mathbb{Z}$; for otherwise the cycle $\mathcal{O}$ misses $(K^-, C^+) + \mathbb{Z}$ and hence $\mathcal{O}$ is a twist lifted $\frac{p}{q}$-cycle for $H_{j;l}$ also and thus $Rot(H_j) = \{\frac{p}{q}\}$ contradicting $\Psi_1 > \Psi_2$. Thus, let $N \in (K^-, C)$ be such that $K^-$ and $N$ are consecutive points in $\mathcal{O}$. Then, (use Fact 5 and (7))

$$H_{j;u}^{jq+q_2}(N) = K^- + jp + p_2 = H_{j;u}^{jq+q_2}(C^+).$$

Therefore $H_{j;u}^{jq+q_2}$ is constant on $[N, C^+]$, contradicting $\mathcal{O}$ missing $[C, C^+] + \mathbb{Z}$. We now have:

$$\Psi_1 > \Psi_2 \Rightarrow \mathsf{L}_{j+1} \neq \mathsf{L}_j \text{ for each } j \geq 0.$$

This completes the proof of Theorem 3. $\square$

*Proof of Theorem 4.* With the notations of Theorems 3 and 4, assume that there exists an $a$ such that $F_{a,b_0}$ is $\frac{p'}{q'}$-locked. To prove Theorem 4, it suffices to consider the cases when $\frac{p'}{q'}$ is a child of $\frac{p}{q}$. We do the case where $\frac{p'}{q'}$ is the right child of $\frac{p}{q}$ as the other case is proven in a similar manner with the $\mathsf{R}_j$'s replacing the $\mathsf{L}_j$'s. Hence, we assume $\frac{p}{q} = \frac{p_1}{q_1} \oplus \frac{p_2}{q_2}$ and $\frac{p'}{q'} = \frac{p}{q} \oplus \frac{p_2}{q_2}$. As we will apply Theorem 3 for both $\frac{p}{q}$ and $\frac{p'}{q'}$ we alter slightly the notation of Theorem 3 to distinguish these applications (for example we will have $\Psi_{i,\frac{p}{q}}$ and $\mathsf{L}_{0,\frac{p}{q}}$ when applying Theorem 3 for $\frac{p}{q}$ and $\Psi_{i,\frac{p'}{q'}}$ and $\mathsf{L}_{0,\frac{p'}{q'}}$ when applying Theorem 3 for $\frac{p'}{q'}$). First notice that using the definition of the $\mathsf{L}_j$'s given in the statement of Theorem 3 and Definition 1 we have that

(29) $$\mathsf{L}_{0,\frac{p}{q}} = \mathsf{L}_{0,\frac{p'}{q'}},$$

since $\mathsf{L}_{0,\frac{p}{q}}$ and $\mathsf{L}_{0,\frac{p'}{q'}}$ are both the intersection of $\mathcal{L}_{\frac{p_2}{q_2}}$ (defined in Definition 2) with the line $b = b_0$ in parameter space. From $\frac{p}{q} < \frac{p'}{q'}$, $F_{a,b_0}$ being $\frac{p'}{q'}$-locked, and Theorem 3 (applied with $\frac{p'}{q'}$) we have that

(30) $$\Psi_{2,\frac{p}{q}} < \Psi_{1,\frac{p'}{q'}} \leq \Psi_{2,\frac{p'}{q'}} < \mathsf{L}_{0,\frac{p'}{q'}},$$

use Remark 3 and the definition of $\Psi_i$ given in Fact 1 for the first inequality in (30), Fact 1(vii) for the second inequality, and Theorem 3 for the third inequality. Combining (29) and (30) we have that $\Psi_{2,\frac{p}{q}} < \mathsf{L}_{0,\frac{p}{q}}$, and therefore (apply Theorem 3 with $\frac{p}{q}$) $\Psi_{1,\frac{p}{q}} < \Psi_{2,\frac{p}{q}}$. Hence, from Fact 1(vii), the line $b = b_0$ intersects the $\frac{p}{q}$-locking region in a segment of positive length. $\square$

*Proof of Theorem 5.* With the notation of Theorem 5 we show that the $T^I_{r_j}$'s converge to a point which lies on the boundary of $L^I_{\frac{p}{q}}$ and which is bounded away from the critical line. The argument for the case of the $T^I_{l_j}$'s is similar. That the $T^I_{r_j}$'s converge to a point, say $A$, follows from Corollary 3. Facts 2 and 5 imply that $A \in Lock_{\frac{p}{q}}$. Moreover, since the map is locked at each $T^I_{r_j}$, and $r_j \neq \frac{p}{q}$ for all $j$, it follows that $A$ lies on the joint boundary of $Lock_{\frac{p}{q}}$ and the $\frac{p}{q}$ Arnold tongue. It remains to show that $A$ is bounded away from the critical line. From Theorem 2 we have that each $T^I_{r_j}$ belongs to $\mathcal{R}_{\frac{p}{q}}$. Hence, $A \in \mathcal{R}_{\frac{p}{q}}$. Suppose $A$ is on the critical line. Since on the critical line $K = C$, it then follows from $A \in \mathcal{R}_{\frac{p}{q}}$ and $F'_{u,v}(K = C) = 0$ that $K = C$ belongs to a stable $\frac{p}{q}$ orbit. This contradicts Remark 3 and Fact 9. Hence we have the first part of Theorem 5.

Let $\{T^I_{\frac{c_j}{d_j}}\}$ be as in the statement of Theorem 5. The sequence of points $\{T^I_{\frac{c_j}{d_j}}\}$ has an accumulation point, $P$, by compactness. By continuity of the rotation interval (Fact 2), the rotation set at P is $\{\omega\}$, and P is on the critical line [B, p. 368] (no $T^I_{\frac{c_j}{d_j}}$ is *below* the critical line). Suppose $P$ is not the limit of the sequence $\{T^I_{\frac{c_j}{d_j}}\}$. Then the sequence $\{T^I_{\frac{c_j}{d_j}}\}$ has an accumulation point, $Q \neq P$. There also, the rotation set is $\{\omega\}$. But there is only one point on the critical line where the rotation set is $\omega$, so we have a contradiction. Thus the point $P$ on the critical line is the limit of the sequence $\{T^I_{\frac{c_j}{d_j}}\}$. $\square$

*Proof of Corollary 2.* With $A$ as in the proof of Theorem 5, consider any horizontal line in the parameter plane from the one containing A down to, but not including, the critical line. By the same argument as given above for tips, the point, $A_b$, of intersection of the horizontal line with the right boundary of the $\frac{p}{q}$ Arnold tongue is the limit of the sequence of intersections of the left (or right) boundaries of the locking regions $T^I_{r_j}$ with the horizontal line. Thus every sufficiently small neighborhood of any of the points, $A_b$, intersects the collection of locking-region boundaries in a set with infinitely many connected components. □

## References

[ALM] L. Alsedá, J. Llibre, and M. Misiurewicz, *Combinatorial dynamics and entropy in dimension one*, Advanced series in nonlinear dynamics. **5**, World Scientific Publishing Co., River Edge NJ, 1993.

[B] P. L. Boyland, *Bifurcations of circle maps: Arnol'd tongues, bistability and rotation intervals*, Commun. Math. Phys. **106** (1986), 353–381.

[BJ] P. M. Blecher and M. V. Jakobson, *Absolutely continuous invariant measures for some maps of the circle*, Statistical Physics and Dynamical Systems; Progress in Physics **10**, Birkhäuser, Boston-Basel-Stuttgart, 1985, pp. 303-315.

[BT] K. M. Brucks and C. Tresser, *A Farey tree organization of locking regions for simple circle maps*, To appear: Proc. Amer. Math. Soc. (1996).

[CGT] A. Chencinier, J. M. Gambaudo, and C. Tresser, *Une remarque sur la structure des endomorphismes de degre une du cercle* **t. 299 serie I** (1984), 771–773.

[EKT] A. Epstein, L. Keen, and C. Tresser, *The set of maps $f_{a,b} : \theta \mapsto \theta + a + \frac{b}{2\pi} \cdot sin(2\pi\theta)$ with any given rotation numbers is contractible*, Commun. Math. Phys. **173** (1995), 313–333.

[FT] B. Friedman and C. Tresser, *Comb structure in hairy boundaries: some transition problems for circle maps*, Phys. Lett. A **117** (1986), 15–22.

[GT] L. Goldberg and C. Tresser, *Rotation orbits and the farey tree*, preprint I. B. M. (1991).

[HW] G. H. Hardy and E. M. Wright, *An introduction to the theory of numbers*, Clarendon Press, Oxford, 1979.

[LT] J. C. Lagarias and C. P. Tresser, *A walk along the branches of the extended Farey tree*, IBM J. Res. Dev. **39** (1995), 283–294.

[MT] R. S. Mackay and C. Tresser, *Transition to topological chaos for circle maps*, Physica D **19** (1986), 206–237.

[P] H. Poincaré, *Sur les courbes définies par des équations différentielles*, J. Math. Pures et Appl. $4^{\grave{e}me}$ série **1** (1885), 167–244, also in "Euvres Complètes, t.1" Gauthier-Villars, Paris (1951).

[RS1] J. Ringland and M. Schell, *The Farey tree embodied — in bimodal maps of interval*, Phys. Letters A **136** (1989), 379–386.

[RS2] J. Ringland, N. Issa and M. Schell, *From U sequence to Farey sequence: a unification of one-parameter scenarios*, Phys. Rev. A **41** (1990), 4223–4235.

[RS3] J. Ringland and M. Schell, *Geneology and bifurcation skeleton for cycles of the interated two-extremum map of the interval*, SIAM J. Math. Analysis **22** (1991), 1354–1371.

[RT] J. Ringland and C.Tresser, *A geneology for the finite kneading sequences of bimodal maps on the interval*, Trans. AMS **347** (1995), 4599–4624.

[STZ] R. M. Siegel, C. Tresser, and G. Zettler, *A decoding problem in dynamics and in number theory*, Chaos **2, no 4** (1982), 473–494.

Department of Mathematical Sciences, U-WI-Milwaukee, Milwaukee, WI 53211
*E-mail address*: kmbrucks@math.uwm.edu

Department of Mathematics, SUNY Buffalo, Buffalo, NY 14214-3093
*E-mail address*: ringland@acsu.buffalo.edu

Thomas J. Watson Research Center, I.B.M., P.O. Box 218, Yorktown Heights, NY 10598
*E-mail address*: tresser@watson.ibm.com